\documentclass[10pt,a4paper]{amsart}
\usepackage[latin1]{inputenc}
\usepackage[T1]{fontenc}
\usepackage[english]{babel}

\usepackage{amsmath,amssymb,amsthm,amsfonts}

\theoremstyle{plain}
\newtheorem{theorem}{Theorem}[section]
\newtheorem{corollary}[theorem]{Corollary}
\newtheorem{lemma}[theorem]{Lemma}
\newtheorem{prop}[theorem]{Proposition}

\theoremstyle{definition}
\newtheorem{remark}[theorem]{Remark}

\newtheorem{problem}[theorem]{Problem}

\newcommand{\C}{\mathbb{C}}

\newcommand{\R}{\mathbb{R}}

\newcommand{\N}{\mathbb{N}}

\newcommand{\T}{\mathbb{T}}
\newcommand{\eps}{\varepsilon}

\DeclareMathOperator{\e}{e}

\DeclareMathOperator{\sign}{sign}
\DeclareMathOperator{\supp}{supp}

\renewcommand{\leq}{\leqslant}
\renewcommand{\geq}{\geqslant}

\newcounter{equi1}

\begin{document}
\title{On the numerical radius of operators in Lebesgue spaces}

\author[Mart\'{\i}n]{Miguel Mart\'{\i}n}
\author[Mer\'{\i}]{Javier Mer\'{\i}}

\address[Mart\'{\i}n \& Mer\'{\i}]{Departamento de An\'{a}lisis Matematico\\
Facultad de Ciencias\\
Universidad de Granada\\
18071 Granada, Spain}

\email{\texttt{mmartins@ugr.es} \qquad \texttt{jmeri@ugr.es}}

\urladdr{http://www.ugr.es/local/mmartins} \urladdr{http://www.ugr.es/local/jmeri}

\author[Popov]{Mikhail Popov}

\address[Popov]{Department of Mathematics\\
Chernivtsi National University\\ str. Kotsjubyn'skogo 2,
Chernivtsi, 58012 Ukraine}

\email{\texttt{misham.popov@gmail.com}}

\curraddr{Department of Mathematics\\ Miami University\\ Oxford, OH\\ 45056 USA E-mail: \texttt{popovm@muohio.edu} }

\urladdr{http://testuvannya.com.ua/mpopov}

\date{November 4th, 2010}

\thanks{First and second authors partially supported by Spanish
MICINN and FEDER project no.\ MTM2009-07498 and Junta de Andaluc\'{\i}a
and FEDER grants G09-FQM-185 and P09-FQM-4911.
Third author supported by Junta de Andaluc\'{\i}a and FEDER grant
P06-FQM-01438 and by Ukr.\ Derzh.\ Tema N 0103Y001103.}


\begin{abstract}
We show that the absolute numerical index of the space $L_p(\mu)$ is $p^{-\frac{1}{p}} q^{-\frac{1}{q}}$ (where $1/p+1/q=1$). In other words, we prove that
$$
\sup\left\{\int |x|^{p-1}|Tx|\, d\mu \, : \ x\in L_p(\mu),\,\|x\|_p=1\right\}
\,\geq \,p^{-\frac{1}{p}} q^{-\frac{1}{q}}\,\|T\|
$$
for every $T\in \mathcal{L}(L_p(\mu))$ and that this inequality
is the best possible when the dimension of $L_p(\mu)$ is
greater than one. We also give lower bounds for the best
constant of equivalence between the numerical radius and the
operator norm in $L_p(\mu)$ for atomless $\mu$ when restricting
to rank-one operators or narrow operators.
\end{abstract}

\maketitle

\section{Introduction and preliminaries} \label{in}

Let $X$ be a real or complex Banach space. Following the
standard notation, by $B_X, \, S_X, \, X^*$ and $\mathcal L(X)$
we denote the closed unit ball, the unit sphere, the dual
space, and the space of all bounded linear operators on $X$
respectively. We write $\T$ for the unit sphere of the base
field $\R$ or $\C$. The \emph{numerical radius} of an operator
$T \in \mathcal L(X)$ is a semi-norm defined as
$$
v(T) = \sup \Bigl\{ \bigl| x^* (Tx) \bigr|\,:
 \ x \in S_X, \ x^* \in S_{X^*}, \ x^*(x) = 1 \Bigr\},
$$
which is obviously smaller or equal than the operator norm. The
\emph{numerical index} of the space $X$ is the constant
$$
n(X) = \inf \bigl\{ v(T)\,: \ T \in \mathcal L(X), \ \|T\| = 1 \bigr\},
$$
equivalently, $n(X)$ is the maximum of those $k\geq 0$ such
that $k\|T\|\leq v(T)$ for every $T\in \mathcal{L}(X)$. This
notion was introduced and studied in the 1970 paper
\cite{D-Mc-P-W}, see also the monographs \cite{B-D1,B-D2} and
the survey paper \cite{KaMaPa} for background. Obviously, $0
\leq n(X) \leq 1$, $n(X) > 0$ means that the numerical radius
is a norm on $\mathcal L(X)$ equivalent to the operator norm
and $n(X)=1$ if and only if numerical radius and operator norm
coincide. It is also not hard to see that $n(X^*) \leq n(X)$,
being the reversed inequality false in general (see \cite[\S
2]{KaMaPa} for a detailed account). There are lots of spaces
with numerical index~$1$ (among classical ones, for instance,
$L_1(\mu)$ and $C(K)$), and some attractive open problems on
them \cite{KaMaPa}. It is interesting to remark that the
numerical index behaves differently in the real and in the
complex cases. So, for every complex Banach space one has that
$n(X) \geq 1/\e$ (and the inequality is the best possible),
nevertheless, $n(X) = 0$ for some real Banach spaces $X$ as
$\ell_2$ or, more in general, for every Hilbert space of
dimension greater than $1$.

The number of Banach spaces whose numerical index is known is
small (see \cite[\S 1]{KaMaPa} for a recent account) and,
therefore, there are many interesting open problems consisting in
calculating, or at least estimating, the numerical index of concrete
Banach spaces. Among classical spaces, one of the most intriguing
open problems is to calculate $n(L_p(\mu))$ for $1<p<\infty$,
$p\neq 2$. Let us fix the notation and terminology on
$L_p$ spaces. Let $(\Omega, \Sigma, \mu)$ be any measure space and
$1< p<\infty$. We write $L_p(\mu)$ for the real or complex Banach
space of (equivalent classes of) measurable scalar functions $x$ defined on $\Omega$ such
that
$$
\|x\|_p=\left(\int_\Omega|x|^p\,d\mu \right)^{\frac{1}{p}}<\infty.
$$
We use the notation $\ell_p^m$ for the $m$-dimensional
$L_p$-space. We write $q=p/(p-1)$ for the conjugate exponent to
$p$. For any $x \in L_p(\mu)$, we denote
$$
x^\# = \begin{cases} |x|^{p-1} \sign(x) & \text{ in the real case}, \\
|x|^{p-1} \sign(\overline{x}) & \text{ in the complex case}, \end{cases}
$$
which is the unique element in $L_q(\mu)\equiv L_p(\mu)^*$ such that
$$
\|x\|_p^p=\|x^\#\|_q^q \qquad \text{and} \qquad \int_\Omega x\,x^\#\ d\mu =
\|x\|_p\,\|x^\#\|_q=\|x\|_p^p.
$$
Observe that, with this notation, one has
$$
v(T) =\sup\left\{\Bigl|\int_\Omega x^\# T x\,d\mu\Bigr|\, : \ x\in S_{L_p(\mu)}\right\}
$$
for every $T\in \mathcal{L}(L_p(\mu))$. Finally, we consider the
constants
\begin{equation}\label{eq:defMp}
M_p=\max_{t\in[0,1]} \frac{|t^{p-1}-t|}{1+t^p}=\max_{t\geq 1} \frac{|t^{p-1}-t|}{1+t^p}\, ,
\end{equation}
(which is the numerical radius of the operator $T(x,y)=(-y,x)$
defined on the real space $\ell_p^2$, see
\cite[Lemma~2]{MarMer-LP} for instance) and
\begin{equation}\label{eq:defkp}
\kappa_p = \max\limits_{\tau > 0} \frac{\tau^{p-1}}{1 + \tau^p} =
\max\limits_{\lambda \in [0,1]} \lambda^\frac{1}{q} (1 -
\lambda)^\frac{1}{p}
 = \frac{1}{p^{1/p} q^{1/q}}
\end{equation}
(which is the numerical radius of the operator $T(x,y)=(y,0)$
defined on the real or complex space $\ell_p^2$, see
\cite[Lemma~2]{MarMer-LP} for instance).

It has been proved recently that, fixed $p$, all
infinite-dimensional $L_p(\mu)$ spaces have the same numerical
index \cite{Eddari,Eddari-Khamsi,Eddari-Khamsi-Aksoy} (see also
\cite{MarMerPopRad-JMAA} for a different approach) and that
$n(L_p(\mu))> 0$ for $p\neq 2$ in the real case \cite{MMP}. On
the way to state the last result, the authors of \cite{MMP}
introduced the so-called absolute numerical radius of an
operator on a $L_p(\mu)$-space as follows. Given a measure
space $(\Omega, \Sigma, \mu)$, $1<p<\infty$ and $T \in \mathcal
L(L_p(\mu))$, the \textit{absolute numerical radius} of $T$ is
the number
\begin{align*}
|v|(T) & = \sup \left\{ |x^*|\bigl(|Tx|\bigr)\, d\mu \, : \  x \in S_{L_p(\mu)}, \
x^* \in S_{L_p(\mu)^*},  \ x^*(x) = 1 \right\} \\ & =
\sup \left\{ \int_\Omega \bigl|x^\# Tx \bigr|\, d\mu \, : \  x \in
S_{L_p(\mu)} \right\} =
\sup \left\{ \int_\Omega |x|^{p-1}|Tx|\,d\mu \, : \  x \in S_{L_p(\mu)}\right\}.
\end{align*}
It is clear that $|v|$ is a seminorm on $\mathcal{L}(L_p(\mu))$
satisfying
$$
v(T)\leq |v|(T) \leq \|T\| \qquad \bigl(T\in \mathcal{L}(L_p(\mu)) \bigr).
$$
In \cite{MMP} it is shown that $n(L_p(\mu))$ is positive by proving that both inequalities above can be reversed up to a positive constant. Namely, it is
shown that
$$
\frac{1}{2\e}\|T\|\leq |v|(T)\qquad \text{ and } \qquad
\frac{M_p}{6}|v|(T)\leq v(T)
$$
for every $T\in \mathcal{L}(L_p(\mu))$, giving $n(L_p(\mu))\geq
\frac{M_p}{12\e}>0$.

For notational convenience, we introduce the definition of the
\emph{absolute numerical index} of $L_p(\mu)$ as the number
\begin{align*}
|n|(L_p(\mu)) & = \inf \bigl\{ |v|(T)\,: \ T \in \mathcal
L(L_p(\mu)), \ \|T\| = 1 \bigr\} \\ & = \max \bigl\{ k\geq 0 \,: \
k\|T\|\leq |v|(T) \ \forall T \in \mathcal{L}(L_p(\mu))\bigr\}
\end{align*}
and the aforementioned result of \cite{MMP} just says that
$|n|(L_p(\mu))\geq \frac{1}{2\e}$. Our first goal in this paper
is to calculate the exact value of $|n|(L_p(\mu))$, namely,
$|n|(L_p(\mu)) = \kappa_p$ (if the dimension of
$L_p(\mu)$ is greater than one) in both the real and the
complex case. In other words, we will prove that,
$$
\sup\left\{\int |x|^{p-1}|Tx|\, d\mu \, : \ x\in L_p(\mu),\,\|x\|_p=1\right\}
\geq \kappa_p\,\|T\|
$$
for every $T\in \mathcal{L}(L_p(\mu))$ and that this inequality
is the best possible when the dimension of $L_p(\mu)$ is
greater than one. As a corollary, we get an improvement of the
estimation of $n(L_p(\mu))$ obtained in \cite{MMP}. Namely, in
the real case, we get
$$
n(L_p(\mu))\geq \frac{\kappa_p\,M_p}{6}.
$$
In other words, in the real case,
$$
\sup\left\{\Bigl|\int |x|^{p-1}\sign(x)\,T x\ d\mu \Bigr| \, : \ x\in L_p(\mu),\,\|x\|_p=1\right\}
\geq \frac{\kappa_p\,M_p}{6}\,\|T\|
$$
for every $T\in \mathcal{L}(L_p(\mu))$.

Next, we study numerical radius of rank-one operators on
$L_p(\mu)$. For notational convenience again, we define the
\emph{rank-one numerical index} of an arbitrary Banach space
$X$ as the number
\begin{align*}
n_1(X) &= \inf \bigl\{ v(T)\,: \ T \in \mathcal L(X), \ \|T\| = 1,\ \text{$T$ rank-one} \bigr\} \\ &=
\max \{ k\geq 0\,:\ k\|T\|\leq v(T)\ \forall T\in \mathcal{L}(X)\ \text{rank-one}\}.
\end{align*}
Our results state that for every atomless measure
$\mu$,
$$
n_1(L_p(\mu))\geq \kappa_p^2
$$
in both the real and the complex cases. This result is not
sharp for values of $p$ close to $2$ as, for instance,
$n_1(L_2(\mu))=\frac12$ if the dimension of $L_2(\mu)$ is
greater than $1$. On the other hand, the estimation for $n_1(L_p(\mu))$
tends to $1$ as $p \to 1$ or $p \to \infty$.

Finally, the last part of the paper is devoted to study numerical
radius of the so-called narrow operators on $L_p(\mu)$ when the
measure $\mu$ is atomless and finite (a class of operators containing compact
operators, see section~\ref{sec:Lpnar} for the definition and
background). Defining the \emph{narrow numerical index} of
$L_p(\mu)$ as
\begin{align*}
n_{\rm nar}(L_p(\mu)) & = \inf \bigl\{ v(T)\,: \ T \in \mathcal L(L_p(\mu)),\
\|T\| = 1, \ T \ \text{narrow} \bigr\} \\
&= \max \{ k\geq 0\,:\ k\|T\|\leq v(T)\ \forall T\in \mathcal{L}(L_p(\mu))\ \text{narrow}\},
\end{align*}
we prove that
$$
n_{\rm nar}(L_p(\mu)) \geq \kappa_p^2 \ \ \text{in the complex case},
\qquad n_{\rm nar}(L_p(\mu)) \geq \max\limits_{\tau > 0} \frac{\kappa_p \tau^{p-1} - \tau}{1 + \tau^p}
\ \ \text{in the real case}.
$$
Notice that the inequality for the real case gives a positive
estimate for $1<p<\infty$ ($p\neq 2$) which tends to $1$ as $p \to 1$ or $p
\to \infty$.

The outline of the paper is as follows. Section~\ref{sec:Lpabs} is devoted to show that $|n|(L_p(\mu))=\kappa_p$. The results on rank-one operators appear in section~\ref{sec:Lprank1} and the results on narrow operators are contained in section~\ref{sec:Lpnar}.

We recall some lattice notation which we will use in the paper. We refer the reader to \cite{AB} for abundant information on lattices and positive operators. Let $E$ be a Banach lattice. For any subset $F \subseteq E$ we write $F^+ = \{x \in F\,: \ x \geq 0\}$. For two elements $x,y \in E$, by $x \vee y$ (resp.\
$x\wedge y$) we denote the least upper bound (resp.\ greatest
lower bound) in $E$ of the two-point set $\{x,y\}$, if it exists.
A linear operator $T:E \longrightarrow E$ is called \emph{positive} provided $T(E^+) \subseteq E^+$, or, in other words, if it sends positive elements to positive elements. An
element $y\in E$ is called a \textit{component} of $x \in E$ if $|y| \wedge |x-y| = 0$. In this case we write $y \sqsubseteq
x$. Let $z\in E$. An element $x \in E$ is called a
$z$-\textit{step function} if $x = \sum_{k=1}^m a_k z_k$ for
some components $(z_k)$ of $z$.

Let $(\Omega, \Sigma, \mu)$ be a measure space. On the real space $L_0(\mu)$ of all (equivalence classes of) $\Sigma$-measurable functions, we consider the ordering $x \leq y$ if and only if $x(t) \leq y(t)$ for almost all $t \in \Omega$. For two functions $x,y \in L_0$, $x \vee y$ (resp.\
$x\wedge y$) is equal to the point-wise maximum (resp.\ minimum) of these functions. For any $x \in L_0(\mu)$ and $A \in \Sigma$ we denote $x_A = x\, \textbf{1}_A$ where $\textbf{1}_A$ is the characteristic function of $A$. The expression $A = B \sqcup C$
for sets $A,B,C \in \Sigma$ means that $A = B \cup C$ and $B
\cap C = \emptyset$. If $E$ is a sublattice of $L_0(\mu)$ and
$x, y \in E$ then $y \sqsubseteq x$ if and only if $y = x\,
\textbf{1}_A$ for some $A \in \Sigma$ and a $\textbf{1}$-step
function is just a simple function {and a $z$-step
function is the product of $z$ by a simple function. In particular, if $x, z \in E$ are simple (= finite valued) functions with $z \geq 0$ and $\supp x \subseteq \supp z$ then $x$ is a $z$-step function.

\section{The absolute numerical index of $L_p(\mu)$} \label{sec:Lpabs}
The main aim of this section is to calculate the absolute
numerical index of the $L_p$ spaces, as shown in the following
result.

\begin{theorem} \label{t:Lpabs}
Let $1<p<\infty$ and let $\mu$ be a positive measure such that
$\dim(L_p(\mu))\geq2$. Then,
$$
|n|(L_p(\mu)) = \kappa_p\,.
$$
\end{theorem}

It is immediate to check that for positive operators on
$L_p(\mu)$, the numerical radius and the absolute numerical radius
coincide. Therefore, the following result is a consequence of the
above theorem. We state here its proof since it is simple and
useful to get a better understanding of the proof of
Theorem~\ref{t:Lpabs}.

\begin{prop}\label{prop-Positive}
Let $1<p<\infty$ and $(\Omega,\Sigma,\mu)$ be a measure space. Then, for
every positive operator $T\in \mathcal{L}(L_p(\mu))$ one has
$$
v(T)\geq\kappa_p\,\|T\|\,.
$$
\end{prop}

\begin{proof}
Let $T \in \mathcal{L}(L_p(\mu))$ be positive with $\|T\| = 1$,
fix $\eps > 0$, and take $x \in S_{L_p(\mu)}$ so that $\|Tx\|^p
\geq (1 - \eps)$ and $x\geq0$ (observe that $x$ can be taken
positive because $T|x|\geq|Tx|$ due to the positivity of $T$).
Next, fix any $\tau>0$, set
$$
y=x\vee\tau Tx \qquad \text{and}\qquad A = \{\omega\in \Omega\, :\, x(\omega) \geq \tau (Tx)(\omega) \},
$$
and observe that
$$
\|y\|^p = \int_A x^p \,d \mu + \int_{\Omega \setminus A} (\tau Tx)^p \,d \mu \leq 1 + \tau^p\qquad \text{and}\qquad y^\#=x^{p-1}\vee(\tau Tx)^{p-1}.
$$
This, together with the positivity of $T$, allows us to write
\begin{align*}
v(T) &\geq \frac{1}{\|y\|^p}\int_{\Omega} y^\# Ty \,d \mu \geq \frac{1}{1 + \tau^p} \int_{\Omega} y^\# Ty \,d \mu \\
&\geq\frac{1}{1 + \tau^p} \int_{\Omega} (\tau Tx)^{p-1} Tx \,d \mu = \frac{\tau^{p-1}}{1 + \tau^p} \int_{\Omega} (Tx)^p  \,d \mu \geq \frac{\tau^{p-1}}{1 + \tau^p} (1 - \varepsilon)
\end{align*}
for every $\tau>0$. Taking supremum on $\tau>0$ and $\eps>0$, we deduce that $v(T)\geq \kappa_p$, as desired.
\end{proof}

The proof of Theorem~\ref{t:Lpabs} depends on the base scalar field. In
the real case it needs some auxiliary results which we state
here. They carry the main idea for the best possible estimation
of the absolute numerical radius in the real case and allow us
to apply positivity arguments to any operator as it has been done in the proof of Proposition~\ref{prop-Positive}.

\begin{lemma} \label{lemma:Lpabs}
Let $E$ be a vector lattice, $z \in E^+$, and $x \in E$ a $z$-step
function with $|x| \leq z$. Then there exist $n\in\N$,
$\lambda_j\in [0,1]$, and $y_j \in E$ with $|y_j|=z$ for
$j=1,\dots,n$ such that $\sum_{j=1}^n \lambda_j=1$ and
$$
x = \lambda_1 y_1 + \cdots + \lambda_n y_n.
$$
\end{lemma}

\begin{proof}
Let $x = \sum_{k=1}^m a_k z_k$ with $a_k \in \R$ and $z_k \sqsubseteq z$, and use induction on $m$. Observe
that the hypothesis $|x|\leq z$ implies that $|a_k|\leq 1$ for every
$k=1,\ldots,m$. For $m=1$, one trivially has that
$x=\frac{1+a_1}{2}z_1+\frac{1-a_1}{2}(-z_1)$. For the induction step assume that the assertion is true for a given $m \in \mathbb N$ and suppose that $x = \sum_{k=1}^{m+1} a_k z_k$ where $z_k \sqsubseteq
z$ and $|a_k| \leq 1$ for $k = 1, \ldots, m+1$. Then for $\widetilde{x} =
\sum_{k=1}^m a_k z_k$ and $\widetilde{z} = z - z_{m+1}\in E^+$ we have that
$z_k \sqsubseteq \widetilde{z}$ for $k = 1, \ldots, m$. By the induction assumption there are $n_0\in\N$, $\widetilde{\lambda}_j\in [0,1]$, and $\widetilde{y}_j \in E$
with $|\widetilde{y}_j| = \widetilde{z}$ for $j=1,\ldots,n_0$ such that $\sum_{j=1}^{n_0}
\widetilde{\lambda}_j = 1$  and $\widetilde{x} = \widetilde{\lambda}_1 \widetilde{y}_1 + \cdots + \widetilde{\lambda}_{n_0} \widetilde{y}_{n_0}$.
Then set $\lambda = \frac{1+a_{m+1}}{2}$ and observe that
\begin{align*}
x&=\widetilde{x}+a_{m+1}z_{m+1}=\lambda(\widetilde{x}+z_{m+1})+(1-\lambda)(\widetilde{x}-z_{m+1})\\
&=\lambda(\widetilde{\lambda}_1 \widetilde{y}_1 +\cdots + \widetilde{\lambda}_{n_0} \widetilde{y}_{n_0} + z_{m+1}) + (1 - \lambda)(\widetilde{\lambda}_1 \widetilde{y}_1 + \cdots + \widetilde{\lambda}_{n_0} \widetilde{y}_{n_0} - z_{m+1})\\
&= \lambda \Bigl(\widetilde{\lambda}_1 (\widetilde{y}_1 + z_{m+1}) + \cdots + \widetilde{\lambda}_{n_0} (\widetilde{y}_{n_0} + z_{m+1}) \Bigr)
+(1 - \lambda)\Bigl(\widetilde{\lambda}_1 (\widetilde{y}_1 - z_{m+1}) + \cdots + \lambda_{n_0} (\widetilde{y}_{n_0} - z_{m+1}) \Bigr).
\end{align*}
Finally, take $n=2n_0$ and
\begin{align*}
\lambda_j&=\lambda\widetilde{\lambda}_j, & y_j&=\widetilde{y}_j+z_{m+1} & \text{for}\qquad  j&=1,\dots,n_0 \qquad \qquad \text{and}\\
\lambda_j&=(1-\lambda)\widetilde{\lambda}_j, & y_j&=\widetilde{y}_j-z_{m+1} & \text{for}\qquad  j&=n_0+1,\dots,2n_0
\end{align*}
which fulfill the desired conditions.
\end{proof}

\begin{corollary} \label{cor:Lpabs1}
Let $E$ be a vector lattice, $f$ a positive linear functional on
$E$, $T: E \longrightarrow E$ a linear operator, $z \in E^+$, and $x \in E$ a
$z$-step function with $|x| \leq z$. Then, there exists $y \in E$
satisfying $|y| = z$ and
$f \bigl( |Ty| \bigr) \geq f \bigl( |Tx| \bigr)$.
\end{corollary}

\begin{proof}
By Lemma~\ref{lemma:Lpabs} there are $n\in\N$, $\lambda_j\in [0,1]$,
and $y_j \in E$ with $|y_j|=z$ for $j=1,\dots,n$ such that
$\sum_{j=1}^n \lambda_j=1$ and $x = \lambda_1 y_1 + \cdots +
\lambda_n y_n$. Then we can write
$$
f \bigl( |Tx| \bigr) \leq f \bigl( \lambda_1 |Ty_1| + \cdots + \lambda_n |T y_n| \bigr) = \lambda_1 f \bigl( |Ty_1| \bigr) + \cdots + \lambda_n f \bigl( |Ty_n| \bigr)
$$
and so, $f \bigl( |Ty_j| \bigr)
\geq f \bigl( |Tx| \bigr)$ for some $j$.
\end{proof}

\begin{corollary} \label{cor:Lpabs2}
Let $E$ be a sublattice of $L_0(\mu)$ for some measure space
$(\Omega,\Sigma,\mu)$ in which the set of all simple functions
is dense, $f \in (E^*)^+$, $T \in \mathcal L(E)$, $\varepsilon
> 0$, $z \in E^+$ and $x \in E$ such that $|x| \leq z$. Then
there exists $y \in E$ satisfying $|y| = z$ and $f \bigl( |Ty|
\bigr) \geq f \bigl( |Tx| \bigr) - \varepsilon$.
\end{corollary}

\begin{proof}
It follows immediately from Corollary~\ref{cor:Lpabs1} and the
continuity of $f$, $| \cdot |$, and $T$.
\end{proof}

We are ready to prove the main result.

\begin{proof}[Proof of Theorem~\ref{t:Lpabs}]
To prove that $|n|(L_p(\mu))\leq p^{-1/p}q^{-1/q}$, it suffices to
construct a norm one operator $T_0\in\mathcal{L}(L_p(\mu))$ with
$|v|(T_0)\leq p^{-1/p}q^{-1/q}$. Indeed, we pick disjoint sets
$A,B \in \Sigma$ with $0<\mu(A),\mu(B)<\infty$ (this is possible
since $\dim(L_p(\mu))\geq 2$) and define $T_0 \in \mathcal
L\bigl(L_p(\mu)\bigr)$ by
\begin{equation}\label{eq:operator-absolute-positive-rankone}
T_0 x = \mu(A)^{-1/q} \mu(B)^{-1/p} \left( \int_A x \, d \mu \right) \textbf{1}_B \qquad \bigl(x\in L_p(\mu)\bigr).
\end{equation}
It is easy to check that $\|T_0\| = 1$. Now we show that $|v|(T_0)
\leq \kappa_p$. Given any $x \in S_{L_p(\mu)}$, we set
$$
\lambda = \|x_B\|^p = \int_B |x|^p \,d\mu
$$
and observe that
$$
\|x_A\|^p = \int_A |x|^p \,d\mu \leq 1 - \lambda.
$$
Thus,
\begin{align*}
\int_\Omega |x|^{p-1} |T_0 x| \,d\mu &
= \int_\Omega |x|^{p-1} \textbf{1}_B |T_0 x| \,d\mu \leq \left( \int_B |x|^{(p-1)q} \,d\mu \right)^{1/q} \left( \int_\Omega |T_0 x|^p \,d\mu \right)^{1/p} \\
& = \left( \int_B |x|^p \,d\mu \right)^{1/q} \left( \mu(A)^{-p/q} \mu(B)^{-1} \Bigl| \int_A x \, \,d\mu \Bigr|^p \mu(B) \right)^{1/p} \\
& \leq \lambda^{1/q} \left( \mu(A)^{-p/q} \Bigl( \int_A |x| \, \,d\mu \Bigr)^p \right)^{1/p} \\
&\leq \lambda^{1/q}  \left( \mu(A)^{-p/q} \mu(A)^{p/q} \|x_A\|^p \right)^{1/p}
\leq \lambda^{1/q} (1 - \lambda)^{1/p} \leq \kappa_p.
\end{align*}
Now, we take supremum with $x\in S_{L_p(\mu)}$ to get
$|v|(T_0)\leq \kappa_p$ as desired.

For the more interesting converse inequality, fix $T \in
\mathcal L(L_p(\mu))$ with $\|T\|=1$, $\varepsilon>0$, and
$\tau > 0$, choose $x\in S_{L_p(\mu)}$ so that $\|Tx\|^p_p \geq
1 - \varepsilon$, and set
$$
A = \{\omega\in \Omega \ : \ |x(\omega)| \geq \tau
|(Tx)(\omega)|\}\qquad  \text{and} \qquad B = \Omega \setminus A.
$$
We split the rest of the proof depending on the base scalar field.

\noindent\textit{$\bullet$ Real case.} Using Corollary
\ref{cor:Lpabs2} for $x$, $z = |x|_A + \tau |Tx|_B$, and
$\displaystyle f(u)=\int_\Omega |Tx|^{p-1} u \,d \mu$ ($u\in L_p(\mu)$), choose $y\in L_p(\mu)$
satisfying $|y| = z$ and $f \bigl( |Ty| \bigr) \geq f \bigl(
|Tx| \bigr) - \varepsilon$. Then
$$
\|y\|^p = \|z\|^p \leq 1 + \tau^p \qquad \text{and} \qquad |y| \geq \tau |Tx|,
$$
and therefore, we can write
\begin{align*}
|v|(T) &\geq\int_\Omega \left|\frac{y^{\#}}{\|y^{\#}\|} T\left(\frac{y}{\|y\|}\right)\right| d \mu =\frac{1}{\|y\|^p} \int_\Omega |y|^{p-1} |Ty|\, d \mu \geq \frac{\tau^{p-1}}{1 + \tau^p} \int_\Omega |Tx|^{p-1} |Ty|\, d \mu\\
&\geq \frac{\tau^{p-1}}{1 + \tau^p}\left( \int_\Omega |Tx|^{p-1} |Tx|\, d \mu -\eps\right)\geq \frac{\tau^{p-1}}{1 + \tau^p} (1 - 2\varepsilon)
\end{align*}
for every $\tau>0$. Finally, the arbitrariness of $\varepsilon > 0$
gives $|v|(T)\geq\max\limits_{\tau > 0} \frac{\tau^{p-1}}{1 + \tau^p}$
and so $|n|(L_p(\mu)) \geq \kappa_p$.

\vspace{0,3 cm}

\noindent\textit{$\bullet$ Complex case.} Since $|x| < \tau
|Tx|$ on $B$, it is possible to find measurable functions
$\theta_1, \theta_2: B \longrightarrow \C$ such that
$$
x(\omega) = \frac{1}{2} \theta_1(\omega) + \frac{1}{2} \theta_2(\omega) \qquad
\text{and} \qquad |\theta_j(\omega)|=\tau|(Tx)(\omega)|
\qquad \bigl(\omega\in B,\ j=1,2 \bigr).
$$
Indeed, for $\omega\in B$ define
\begin{align*}
\theta_1(\omega)&=\sign
\big(x(\omega)\big)\left(|x(\omega)|+i\left(\tau^2|(Tx)(\omega)|^2-|x(\omega)|^2\right)^{1/2}\right)\\
\theta_2(\omega)&=\sign
\big(x(\omega)\big)\left(|x(\omega)|-i\left(\tau^2|(Tx)(\omega)|^2-|x(\omega)|^2\right)^{1/2}\right)
\end{align*}
if $x(\omega)\neq0$ and $\theta_1(\omega)=1$,
$\theta_2(\omega)=-1$ if $x(\omega)=0$. Then define
$$
y_j = x_A + \widetilde{\theta}_j \quad (j = 1,2)
$$
where $\widetilde{\theta}_j=\theta_j$ on $B$ and
$\widetilde{\theta}_j=0$ on $A$, and observe that
$$
x=\frac{1}{2} y_1 + \frac{1}{2} y_2, \qquad  \|y_j\|^p \leq 1 +
\tau^p, \qquad \text{and} \qquad  |y_j|= |x|_A + |\widetilde{\theta}_j| \geq
\tau |Tx|.
$$
Therefore, we can write
\begin{align*}
|v|(T) &\geq \frac{1}{2} \, \frac{1}{\|y_1\|^p} \int_\Omega |y_1|^{p-1} |Ty_1| d \mu + \frac{1}{2} \, \frac{1}{\|y_2\|^p} \int_\Omega |y_2|^{p-1} |Ty_2| d \mu\\
&\geq \frac{\tau^{p-1}}{1 + \tau^p} \int_\Omega |Tx|^{p-1} \Bigl( \frac{1}{2} |Ty_1| + \frac{1}{2} |Ty_2| \Bigr) d \mu \geq \frac{\tau^{p-1}}{1 + \tau^p} \int_\Omega |Tx|^{p-1} \left| T \Bigl(\frac{1}{2} y_1 + \frac{1}{2} y_2 \Bigr) \right| d \mu\\
&=\frac{\tau^{p-1}}{1 + \tau^p} \int_\Omega |Tx|^p d \mu \geq \frac{\tau^{p-1}}{1 + \tau^p} (1 - \varepsilon)
\end{align*}
for every $\tau>0$. The arbitrariness of $\varepsilon$ gives us that
$|v|(T)\geq\max\limits_{\tau>0}\frac{\tau^{p-1}}{1 + \tau^p}$ and
hence $|n|(L_p(\mu))\geq \kappa_p$ which finishes the proof.
\end{proof}

We can use Theorem~\ref{t:Lpabs} together with
\cite[Theorem~1]{MMP} to improve the estimation of
$n(L_p(\mu))$ given for the real case in
\cite[Corollary~3]{MMP}.

\begin{corollary}
Let $1<p<\infty$ and $\mu$ a positive mesure. Then, in the real case, one has
$$
n(L_p(\mu))\geq \frac{M_p\kappa_p}{6}\,.
$$
\end{corollary}

\begin{remark}
{\slshape From the proof of Theorem~\ref{t:Lpabs} we deduce that
$\kappa_p$ is also the best constant of equivalence between the
norm and the numerical radius for positive operators (i.e.\ the
inequality in Proposition~\ref{prop-Positive} is the best
possible).\,} This is because the operator defined on
Equation~\ref{eq:operator-absolute-positive-rankone} is clearly
positive.
\end{remark}

\section{The numerical radius of rank-one operators on $L_p(\mu)$} \label{sec:Lprank1}
This section is devoted to estimate the numerical radius of
rank-one operators on $L_p(\mu)$ for atomless measures $\mu$.

\begin{theorem} \label{t:Lprankone}
Let $(\Omega,\Sigma,\mu)$ be an atomless measure space. Then
for $1<p<\infty$ one has
$$
\kappa_p\geq  n_1(L_p(\mu)) \geq \kappa_p^2\,.
$$
\end{theorem}

We need the following easy observation.

\begin{remark}\label{rem:lambda-partition}
{\slshape Let $(\Omega,\Sigma,\mu)$ be an atomless measure
space and let $f_1, \ldots, f_n$ be simple functions on
$\Omega$. Then, given any $\lambda \in [0,1]$, there exists a
partition $\Omega = A \sqcup B$ into measurable subsets such
that
$$
\int_A f_j\, d \mu = \lambda \int_{\Omega} f_j\, d \mu \qquad \text{and}
\qquad \int_B f_j\, d \mu = (1 - \lambda) \int_{\Omega} f_j\, d \mu
$$
for every $j = 1, \ldots, n$.\ } To see that this is true, let
$C_0,C_1,\ldots,C_m$ be a partition of $\Omega$ with $0<\mu(C_k)<\infty$ for $k=1,\ldots,m$, such that all the
functions $f_1, \ldots, f_n$ are null on $C_0$ and constant on every $C_k$. Then,
for each $k=1,\ldots,m$, take $A_k, B_k\in \Sigma$ satisfying
$C_k=A_k\sqcup B_k$, $\mu(A_k)=\lambda \mu(C_k)$, and
$\mu(B_k)=(1-\lambda)\mu(C_k)$, and observe that the sets given
by
$$
A=C_0\cup \bigcup_{k=1}^{m} A_k \qquad \text{and} \qquad B=\bigcup_{k=1}^{m} B_k
$$
form the desired partition of $\Omega$.
\end{remark}

\begin{proof}[Proof of Theorem~\ref{t:Lprankone}]
The first inequality follows from the fact that in the proof of
Theorem~\ref{t:Lpabs} it is constructed a positive and rank-one
operator $T_0$ (see \eqref{eq:operator-absolute-positive-rankone})
such that $\|T_0\|=1$ and $v(T_0)=|v|(T_0)\leq \kappa_p$.
Therefore, $n_1(L_p(\mu))\leq \kappa_p$.

We now prove the more interesting second inequality. Let $T \in
\mathcal L(L_p(\mu))$ be a rank-one operator with norm one, that
is,
$$
Tz = \left( \int_{\Omega} x^\# z\, d \mu \right) y
\qquad \ \bigl(z\in L_p(\mu) \bigr)
$$
for some fixed $x,y \in S_{L_p(\mu)}$. Without loss of
generality, we may and do assume that $x,y$ are simple
functions. Fix $\tau > 0$, $\lambda \in [0,1]$ and set
$$
\theta = \mbox{sign} \, \left( \int_{\Omega} x^\# y \,d \mu \right)\in\T
\qquad \text{and} \qquad \delta = \left| \int_{\Omega} x^\# y \,d \mu \right|.
$$
Using Remark~\ref{rem:lambda-partition}, choose a partition
$\Omega = A \sqcup B$ so that
\begin{align}\label{eq:partition1}
\int_A x^\# y \,d \mu &= \lambda \int_{\Omega} x^\# y \,d \mu = \lambda \theta \delta, &\bigl\| x_A \bigr\|^p = \bigl\| y_A \bigr\|^p &= \lambda,\\
\int_B x^\# y \,d \mu &= (1-\lambda) \int_{\Omega} x^\# y \,d \mu = (1-\lambda) \theta \delta,\qquad \text{and}  &\bigl\| x_B \bigr\|^p = \bigl\| y_B \bigr\|^p &= 1-\lambda.\notag
\end{align}
Then define $\displaystyle{z = \lambda^{-\frac{1}{p}} x_A + \overline{\theta}
(1 - \lambda)^{-\frac{1}{p}} \tau y_B}$ and observe that
\begin{equation*}
\|z\|^p =  \lambda^{-1} \bigl\|x_A\bigr\|^p + (1 - \lambda)^{-1} \tau^{p} \bigl\|y_B\bigr\|^p = 1 + \tau^p.
\end{equation*}
Therefore, we can write
\begin{align}\label{thm-rank-one:radius}
v(T) \geq\left|\int_{\Omega}\frac{z^\#}{\|z^\#\|}T\left(\frac{z}{\|z\|}\right)\right|\,d\mu=
 \frac{1}{\|z\|^p} \left| \int_{\Omega} z^\# Tz \, d \mu \right|
 = \frac{1}{1+\tau^p} \left| \int_{\Omega} x^\# z \,d \mu \right| \, \left| \int_{\Omega} z^\# y \,d \mu \right|.
\end{align}
Besides, using the fact that $(u+v)^\# = u^\# + v^\#$ for
disjointly supported elements $u,v \in L_p(\mu)$\,, it is clear that
$z^\# = \lambda^{-\frac{1}{q}} x_A^\# + \theta (1 -
\lambda)^{-\frac{1}{q}} \tau^{p-1} y_B^\#$. Using this and
\eqref{eq:partition1} it is easy to check that
\begin{align*}
\left| \int_{\Omega} x^\# z \,d \mu \right|
&=\left|\lambda^{-\frac{1}{p}} \int_A x^\# x \,d \mu +
\overline{\theta} (1 - \lambda)^{-\frac{1}{p}} \tau \int_B x^\# y\,d\mu\right|\\
&=\left|\lambda^{-\frac{1}{p}} \lambda +
\overline{\theta} (1 - \lambda)^{-\frac{1}{p}} \tau (1 - \lambda) \theta \delta\right| = \lambda^\frac{1}{q} + (1 - \lambda)^\frac{1}{q} \tau \delta \geq \lambda^\frac{1}{q}
\end{align*}
and
\begin{align*}
\left| \int_{\Omega} z^\# y \,d \mu \right| &=
\left| \lambda^{-\frac{1}{q}} \int_A x^\# y \,d \mu +
\theta (1 - \lambda)^{-\frac{1}{q}} \tau^{p-1} \int_B y^\# y \,d \mu \right|\\
&= \left| \lambda^{-\frac{1}{q}} \lambda \theta \delta + \theta (1 - \lambda)^{-\frac{1}{q}} \tau^{p-1} (1 - \lambda) \right|
= \lambda^{\frac{1}{p}} \delta + (1 - \lambda)^\frac{1}{p} \tau^{p-1} \\&\geq (1 - \lambda)^\frac{1}{p} \tau^{p-1}.
\end{align*}
This, together with \eqref{thm-rank-one:radius}, tells us that
$$
v(T)\geq \lambda^\frac{1}{q}(1 - \lambda)^\frac{1}{p}\frac{\tau^{p-1}}{1 + \tau^p}
$$
for every $\tau>0$ and every $\lambda\in[0,1]$. Finally, since
$\max\limits_{\tau>0}\frac{\tau^{p-1}}{1 +
\tau^p}=\max\limits_{\lambda\in[0,1]}\lambda^\frac{1}{q}(1 -
\lambda)^\frac{1}{p}=\kappa_p$, one obtains that $v(T)\geq\kappa_p^2$
which finishes the proof.
\end{proof}

Note that for $p \to 1$ and $p \to \infty$ this gives the best
possible estimation of the order of $n_1(L_p(\mu))$ because
$\kappa_p^2 \to 1$. Nevertheless, for $p = 2$ one has
$\kappa_2^2 = 1/4$, while the best estimation is $1/2$, as the
following easy result shows.

\begin{prop} \label{prop:L2}
Let $H$ be a real Hilbert space of dimension greater than one. Then
$n_1(H) = \frac{1}{2}$.
\end{prop}

\begin{proof}
We fix a rank-one operator $T\in L(H)$ with $\|T\|=1$. Then, $T$ has
the form $Tx=(x\mid x_1)\,x_2$ for some elements $x_1,x_2\in S_H$.
If $|(x_1\mid x_2)|=1$ then $v(T)\geq |(Tx_1\mid x_1)|=1$ and we are done. If
otherwise $|(x_1\mid x_2)|<1$, take $x=\frac{x_1+\theta
x_2}{\|x_1+\theta x_2\|}\in S_H$ for $\theta\in\{-1,1\}$ and observe
that
\begin{align*}
v(T) \geq |(Tx\mid x)| &=\left|\frac{\big(x_1+\theta x_2\mid
x_1\big)\big(x_2\mid x_1+\theta x_2\big)} {\|x_1+\theta x_2\|^2}\right|\\
&=\left| \frac{\big[1 + \theta(x_2\mid x_1)\big]^2}{\|x_1+\theta
x_2\|^2}\right|=\left|\frac{\big[1 + \theta(x_2\mid x_1)\big]^2}{2\big[1 + \theta(x_2\mid x_1)\big]}\right| = \left|\frac{1+\theta (x_2\mid x_1)}{2}\right|.
\end{align*}
By just choosing the suitable $\theta\in\{-1,1\}$ one obtains
$v(T)\geq 1/2$ and so $n_1(H)\geq \frac12$.

For the converse inequality, observe that if we take $x_1,x_2$
orthogonal, then for each $x\in S_H$ one has that $(x\mid
x_1)^2+(x\mid x_2)^2\leq 1$ and, therefore,
$$
|(Tx\mid x)| =|(x\mid x_1)|\,|(x\mid x_2)|=\frac{(x\mid x_1)^2+(x\mid x_2)^2}{2}-\frac12\big(|(x\mid x_1)|-|(x\mid x_2)|\big)^2\leq\frac12
$$
which implies $v(T)\leq \frac12$ and so $n_1(H)\leq \frac12$\,.
\end{proof}

\section{The numerical radius of narrow operators} \label{sec:Lpnar}
In Section~\ref{sec:Lprank1} we obtained an estimate for the
numerical radius of rank-one operators in $L_p(\mu)$, it is natural
to ask if it is possible to obtain a similar estimate for
finite-rank operators. The aim of this section is to prove that
it is so. In fact, we will do the work for the wider class of
narrow operators. Let us recall the relevant definitions. An
operator $T \in \mathcal L (E,X)$ on a (real or complex) K\"{o}the
function space $E$ on a finite measure space $(\Omega, \Sigma,
\mu)$ acting to a Banach space $X$ is \textit{narrow} if for
each $A \in \Sigma$ and each $\varepsilon > 0$ there is an $x
\in E$ such that $x^2 = \textbf{1}_A$, $\displaystyle \int_{\Omega}x \,
d \mu = 0$ and $\|Tx\| < \varepsilon$. The conditions $x^2 =
\textbf{1}_A$, $\displaystyle \int_{\Omega}x \, d \mu = 0$ mean that
there exists a decomposition $A = A^+ \sqcup A^-$ into sets of
equal measure with $x = \textbf{1}_{A^+} - \textbf{1}_{A^-}$. This concept was
introduced in \cite{PP} and developed in some other papers
\cite{KKW,KP03,MMP1} (see also the expository paper \cite{P}).
Note that if $A \in \Sigma$ is an atom then $T\boldsymbol{1}_A
= 0$ for any narrow operator $T \in \mathcal L(E,X)$, thus, the
notion of narrow operator is nontrivial only for atomless
measure spaces $(\Omega, \Sigma, \mu)$. For a more general
consideration of narrow operators we refer the reader to
\cite{MMP1}. If the norm of $E$ is absolutely continuous, then
for every Banach space $X$ every compact ($AM$-compact,
Dunford-Pettis\ldots) operator $T \in \mathcal{L}(E,X)$ is
narrow \cite{MMP1,PP}. For $E = L_p(\mu)$ this is easy to see
using the technique of the Rademacher system. Indeed, consider
any Rademacher system $(r_n)$ on $A$ \cite{PP}. Since $(r_n)$
is a weakly null sequence, we have that $\|Tr_n\| \to 0$ as $n
\to \infty$. However, the converse is not true: there exists a
narrow projection $P \in \mathcal L (L_p[0,1])$ of norm one
onto a subspace of $L_p[0,1]$ isometric to $L_p[0,1]$
\cite{PP}.

Our estimate for the numerical radius of narrow operators in
$L_p(\mu)$ depends on the base scalar field. For the complex
case we obtain the same estimate as we did for rank-one
operators.

\begin{theorem} \label{t:Lpnar}
Let $(\Omega,\Sigma,\mu)$ be an atomless finite measure space.
Then, for every $1<p<\infty$ one has
$$
n_{\rm nar}(L_p(\mu)) \geq \kappa_p^2 \ \ \text{in the complex case},
\qquad n_{\rm nar}(L_p(\mu)) \geq \max\limits_{\tau > 0} \frac{\kappa_p \tau^{p-1} - \tau}{1 + \tau^p}
\ \ \text{in the real case}.
$$
\end{theorem}

Notice that the inequality for the real case gives a positive
estimate for $1<p<\infty$ ($p\neq 2$) which tends to $1$ as $p \to 1$ or $p
\to \infty$.

To prove this result we need the following lemmas which suggest
that a narrow operator behaves almost as a rank-one operator
when it is restricted to a suitable finite dimensional subspace
of arbitrarily large dimension.

\begin{lemma} \label{lem:nar1}
Let $(\Omega,\Sigma,\mu)$ be an atomless finite measure space,
$1 \leq p < \infty$, $T \in \mathcal L(L_p(\mu))$ a narrow
operator, $x\in L_p(\mu)$ a simple function, $Tx = y$, and
$\Omega = D_1 \sqcup \ldots \sqcup D_\ell$ any partition and
$\varepsilon > 0$. Then there exists a partition $\Omega = A
\sqcup B$ such that
\begin{enumerate}
\item[(\emph{i})] $\bigl\| x_A \bigr\|^p =\bigl\| x_B \bigr\|^p= 2^{-1} \|x\|^p$.
\item[(\emph{ii})] $\mu(D_j \cap A) = \mu(D_j \cap B) = \frac{1}{2} \mu(D_j)$ for each $j = 1, \ldots, \ell$.
\item[(\emph{iii})] $\bigl\| Tx_A - 2^{-1} y \bigr\| < \varepsilon$ and $\bigl\| Tx_B - 2^{-1} y \bigr\| < \varepsilon$.
\end{enumerate}
\end{lemma}

\begin{proof}
Let $x = \sum\limits_{k=1}^m a_k \textbf{1}_{C_k}$ for some $a_k \in
\mathbb K$ and $\Omega = C_1 \sqcup \ldots \sqcup C_m$. For each $k
= 1, \ldots, m$ and $j = 1, \ldots, \ell$ define sets $E_{k,j} = C_k
\cap D_j$ and, using the definition of narrow operator, choose
$u_{k,j} \in L_p(\mu)$ so that
$$
u_{k,j}^2 = \textbf{1}_{E_{k,j}}, \qquad \int_{\Omega} u_{k,j} \,d \mu
= 0, \qquad \text{and}\qquad  |a_k|\|Tu_{k,j}\| < \frac{2
\varepsilon}{m\ell}\,.
$$
Then set
$$
E_{k,j}^+ = \bigl\{ t \in E_{k,j} \ : \ u_{k,j}(t) \geq 0 \bigr\},
\qquad E_{k,j}^- = E_{k,j} \setminus E_{k,j}^+
$$
which satisfy $\mu(E_{k,j}^+)=\mu(E_{k,j}^-)=\frac12\mu(E_{k,j})$,
and define
$$
A = \bigcup\limits_{k=1}^m \bigcup\limits_{j=1}^\ell E_{k,j}^+
\qquad \text{and} \qquad B = \bigcup\limits_{k=1}^m
\bigcup\limits_{j=1}^\ell E_{k,j}^-.
$$
Let us show that the partition $\Omega = A \sqcup B$ has the desired
properties.  Indeed, observe that
\begin{align*}
\bigl\| x_A \bigr\|^p &= \sum\limits_{k=1}^m \sum\limits_{j=1}^\ell
|a_k|^p \mu(E_{k,j}^+) = \sum\limits_{k=1}^m |a_k|^p
\sum\limits_{j=1}^\ell \frac{1}{2} \mu(E_{k,j}) = \frac{1}{2}
\sum\limits_{k=1}^m |a_k|^p \mu(C_k) = \frac{1}{2} \|x\|^p
\end{align*}
and that one obviously has $\bigl\| x_B \bigr\|^p=\bigl\| x_A
\bigr\|^p$, thus $(i)$ is proved.

Since $E_{k,j}^+ \subseteq E_{k,j} \subseteq D_j$, for each
$j_0\in\{1,\ldots,\ell\}$ we have that
$$
D_{j_0} \cap A = \bigcup\limits_{k=1}^m \bigcup\limits_{j=1}^\ell
\Bigl(D_{j_0} \cap E_{k,j}^+ \Bigr) = \bigcup\limits_{k=1}^m
E_{k,j_0}^+
$$
and hence
$$
\mu \bigl( D_{j_0} \cap A \bigr) = \sum\limits_{k=1}^m \mu
\bigl(E_{k,j_0}^+ \bigr) = \frac{1}{2} \sum\limits_{k=1}^m \mu
\bigl( E_{k,j_0} \bigr) = \frac{1}{2} \sum\limits_{k=1}^m \mu \bigl(
C_k \cap D_{j_0} \bigr) = \frac{1}{2} \mu(D_{j_0}).
$$
Analogously it is proved that $\mu(D_j\cap B)=\frac12\mu(D_j)$ for
every $j\in\{1,\dots,\ell\}$ which finishes $(ii)$. To prove $(iii)$
observe that
$$
x_A - x_B = \sum\limits_{k=1}^m \sum\limits_{j=1}^\ell a_k \left(
\textbf{1}_{E_{k,j}^+} - \textbf{1}_{E_{k,j}^-} \right) =
\sum\limits_{k=1}^m \sum\limits_{j=1}^\ell a_k u_{k,j}
$$
and hence,
$$
\Bigl\| T (x_A - x_B) \Bigr\| \leq \sum\limits_{k=1}^m
\sum\limits_{j=1}^\ell |a_k|\|T u_{k,j}\| < 2 \varepsilon.
$$
Therefore, one has that
$$
\Bigl\| T x_A - \frac{1}{2} y \Bigr\| = \frac{1}{2} \Bigl\| 2 Tx_A - Tx_A - Tx_B \Bigr\| = \frac{1}{2} \Bigl\| T (x_A- x_B) \Bigr\| < \varepsilon.
$$
Analogously, one obtains $\displaystyle{\Bigl\| Tx_B - \frac{1}{2} y
\Bigr\| < \varepsilon}$ finishing the proof of $(iii)$.
\end{proof}

\begin{lemma} \label{lem:nar2}
Let $(\Omega,\Sigma,\mu)$ be an atomless finite measure space,
$1 \leq p < \infty$, let $T \in \mathcal L(L_p(\mu))$ be a
narrow operator, and let $x, y \in L_p(\mu)$ be simple functions
such that $Tx = y$. Then for each $n \in \mathbb N$ and
$\varepsilon > 0$ there exists a partition $\Omega = A_1 \sqcup
\ldots \sqcup A_{2^n}$ such that for each $k = 1, \ldots, 2^n$
one has
\begin{enumerate}
  \item $\bigl\| x_{A_k} \bigr\|^p = 2^{-n} \|x\|^p$.
  \item $\bigl\| y_{A_k} \bigr\|^p = 2^{-n} \|y\|^p$.
  \item $\bigl\| Tx_{A_k} - 2^{-n} y \bigr\| < \varepsilon$.
\end{enumerate}
\end{lemma}

\begin{proof}
Let $y = \sum\limits_{j=1}^\ell b_j \textbf{1}_{D_j}$ for some
$b_j \in \mathbb K$ and $\Omega = D_1 \sqcup \ldots \sqcup
D_\ell$. We proceed by induction on $n$. Suppose first that $n
= 1$ and use Lemma \ref{lem:nar1} to find a partition $\Omega =
A \sqcup B$ satisfying properties $(i) - (iii)$. Then $(i)$ and
$(iii)$ mean $(1)$ and $(3)$ for $A_1 = A$, $A_2 = B$. Besides,
observe that $(2)$ follows from $(ii)$:
$$
\bigl\| y_{A_1} \bigr\|^p = \sum\limits_{j=1}^\ell |b_j|^p \mu(A_1 \cap D_j) = \sum\limits_{j=1}^\ell |b_j|^p \frac{1}{2} \mu(D_j) = \frac{1}{2} \|y\|^p
$$
and analogously $\bigl\| y_{A_2} \bigr\|^p = 2^{-1} \|y\|^p$.

For the induction step suppose that the statement of the lemma
is true for $n \in \mathbb N$ and find a partition $\Omega =
A_1 \sqcup \ldots \sqcup A_{2^n}$ such that for every
$k=1,\dots,2^n$ the following hold:
\begin{equation}\label{eq:lemma-narrow-2-induction}
\bigl\| x_{A_k} \bigr\|^p = 2^{-n} \|x\|^p, \qquad
\bigl\| y_{A_k} \bigr\|^p = 2^{-n} \|y\|^p, \qquad \text{and}\qquad
\bigl\| Tx_{A_k} - 2^{-n} y \bigr\| < \eps.
\end{equation}
Then, for each $k = 1,
\ldots, 2^n$ use Lemma \ref{lem:nar1} for $x_{A_k}$ instead of $x$,
$Tx_{A_k}$ instead of $y$, the decomposition
$$
\Omega = \bigsqcup\limits_{k=1}^{2^n} \bigsqcup\limits_{j=1}^\ell (D_j \cap A_k)
$$
instead of $\Omega = D_1 \sqcup \ldots \sqcup D_\ell$ and
$\frac{\varepsilon}{2}$ instead of $\varepsilon$, and find a
partition $\Omega = A(k) \sqcup B(k)$ satisfying properties
$(i) - (iii)$. Namely, for each $k=1,\ldots,2^n$ we have that:
\begin{enumerate}
  \item[(\emph{i})] $\bigl\| x_{(A_k \cap A(k))} \bigr\|^p=\bigl\| x_{(A_k \cap B(k))} \bigr\|^p = 2^{-1} \bigl\| x_{A_k} \bigr\|^p$.
  \item[(\emph{ii})] $\mu(D_j \cap A_k \cap A(k)) = \mu(D_j \cap A_k \cap B(k)) = \frac{1}{2} \mu(D_j \cap A_k)$ for each $j = 1, \ldots, \ell$.
  \item[(\emph{iii})] $\bigl\| Tx_{(A_k \cap A(k))} - 2^{-1} Tx_{A_k} \bigr\| < \frac{\varepsilon}{2}$ and $\bigl\| Tx_{(A_k \cap B(k))} - 2^{-1} Tx_{A_k} \bigr\| < \frac{\varepsilon}{2}$.
\end{enumerate}
Let us show that the partition
$$
\Omega = \big(A_1 \cap A(1)\big) \sqcup \ldots \sqcup \big(A_{2^n} \cap A(2^n)\big) \sqcup \big(A_1 \cap B(1)\big) \sqcup \ldots \sqcup \big(A_{2^n} \cap B(2^n)\big)
$$
has the desired properties for $n+1$:

\noindent Property $(1)$: using ($i$) and
\eqref{eq:lemma-narrow-2-induction}, one obtains
$$
\bigl\| x_{(A_k \cap A(k))} \bigr\|^p = \bigl\| x_{(A_k \cap B(k))} \bigr\|^p = 2^{-1} \bigl\| x_{A_k} \bigr\|^p = 2^{-(n+1)} \|x\|^p.
$$

\noindent Property $(2)$: for each $k=1,\ldots,2^n$ use $(ii)$
and \eqref{eq:lemma-narrow-2-induction} to obtain
\begin{align*}
\bigl\| y_{(A_k \cap A(k))} \bigr\|^p = \sum\limits_{j=1}^\ell |b_j|^p \mu\big(D_j \cap A_k \cap A(k)\big) = \frac{1}{2} \sum\limits_{j=1}^\ell |b_j|^p \mu(D_j \cap A_k)
= \frac{1}{2} \bigl\| y_{A_k} \bigr\|^p = 2^{-(n+1)} \|y\|^p
\end{align*}
and analogously $\bigl\| y_{(A_k \cap B)} \bigr\|^p = 2^{-(n+1)}
\|y\|^p$.

\noindent Property $(3)$: for each $k=1,\ldots,2^n$ use $(iii)$
and \eqref{eq:lemma-narrow-2-induction} to write
$$
\Bigl\| Tx_{(A_k \cap A(k))} - 2^{-(n+1)} y \Bigr\| \leq
\Bigl\| Tx_{(A_k \cap A(k))} - 2^{-1} Tx_{A_k} \Bigr\| + \frac{1}{2} \Bigl\| Tx_{A_k} - 2^{-n} y \Bigr\| < \frac{\varepsilon}{2} + \frac{\varepsilon}{2} = \varepsilon
$$
and analogously $\Bigl\| Tx_{(A_k \cap B(k))} - 2^{-(n+1)} y \Bigr\| < \varepsilon$, which completes the proof.
\end{proof}

\begin{lemma} \label{lem:nar3}
Let $(\Omega,\Sigma,\mu)$ be an atomless finite measure space,
$1 \leq p < \infty$, let $T \in \mathcal L(L_p(\mu))$ be a
narrow operator, and let $x, y \in L_p(\mu)$ be simple
functions such that $Tx = y$. Then for each $n \in \mathbb N$,
each number $\lambda$ of the form $\lambda = \frac{j}{2^n}$
where $j \in \{1, \ldots, 2^n-1\}$ and each $\varepsilon > 0$
there exists a partition $\Omega = A \sqcup B$ such that:
\begin{enumerate}
  \item[(A)] $\bigl\| x_A \bigr\|^p = \lambda \|x\|^p$.
  \item[(B)] $\bigl\| y_B \bigr\|^p = (1 - \lambda) \|y\|^p$.
  \item[(C)] $\bigl\| Tx_A - \lambda y \bigr\| < \varepsilon$.
\end{enumerate}
\end{lemma}

\begin{proof}
Use Lemma \ref{lem:nar2} to choose a partition $\Omega = A_1 \sqcup \ldots \sqcup A_{2^n}$ satisfying properties $(1) - (3)$ with $\varepsilon/j$ instead
of $\varepsilon$. Then, setting
$$
A = \bigsqcup\limits_{k=1}^{j} A_k \qquad \text{and} \qquad B = \bigsqcup\limits_{k=j+1}^{2^n} A_k,
$$
one obtains
\begin{align*}
\bigl\| x_A \bigr\|^p &= \sum\limits_{k=1}^{j} \bigl\| x_{A_k} \bigr\|^p = \sum\limits_{k=1}^{j} 2^{-n} \|x\|^p = \lambda \|x\|^p,\\
\bigl\| y_B \bigr\|^p &= \sum\limits_{k=j+1}^{2^n} \bigl\| y_{A_k} \bigr\|^p = \sum\limits_{k=j+1}^{2^n} 2^{-n} \|y\|^p = (1 - \lambda) \|y\|^p,\\
\intertext{and}
\bigl\| Tx_A - \lambda y \bigr\| &= \Bigl\|
\sum\limits_{k=1}^{j} Tx_{A_k} - \sum\limits_{k=1}^{j} 2^{-n}
y \Bigr\| \leq \sum\limits_{k=1}^{j} \bigl\| Tx_{A_k} - 2^{-n}
y \bigr\| < j \frac{\varepsilon}{j} = \varepsilon
\end{align*}
as desired.
\end{proof}

\begin{proof}[Proof of Theorem \ref{t:Lpnar}]
Let $T \in \mathcal L(L_p(\mu))$ be a narrow operator of norm one. Fix
$\varepsilon > 0, \tau > 0,n\in\N$ and $\lambda \in ]0,1[$ of the
form $\lambda = \frac{j}{2^n}$ where $j \in \{1, \ldots, 2^n-1\}$.
Pick a simple function $x \in S_{L_p(\mu)}$ so that $y = Tx$
satisfies $\|y\|^p \geq 1 - \varepsilon$. Without loss of
generality we may assume that $y$ is a simple function since one
can approximate $T$ by a sequence of narrow operators with the
desired property (indeed, take a sequence of simple functions
$(y_m)$ converging to $y$ and define $T_m=T-x^\#\otimes(y-y_m)$.
Then, $T_m(x)=y_m$, $\|T_m-T\|\leq\|y-y_m\|$, and $T_m$ is narrow
for every $m\in\N$ since it is the sum of a rank-one operator and
a narrow one \cite[Proposition~6 on p.~59]{PP}).

Use Lemma~\ref{lem:nar3} to find a partition $\Omega = A \sqcup
B$ satisfying (A)--(C) and use (B) and (C) to obtain the
following estimate:

\begin{align}\label{eq:thm-nar-estimate}
\left| \int_B y^\# T x_{A} \,d \mu - \lambda (1 - \lambda)\|y\|^p
\right| = \left| \int_B y^\# T x_{A} \,d \mu - \lambda \int_B y^\# y\,d\mu
\right|
\leq \bigl\| T x_{A} - \lambda
y \bigr\| < \varepsilon\,.
\end{align}
Then for $\theta\in \T$ define  $z_\theta = \lambda^{-
\frac{1}{p}} x_{A} + \theta (1 - \lambda)^{- \frac{1}{p}} \tau
y_{B}$ and observe, using (A) and (B) of Lemma~\ref{lem:nar3},
that
\begin{equation*}
\|z_\theta\|^p = \lambda^{-1} \bigl\|x_{A} \bigr\|^p + (1 - \lambda)^{-1}
\tau^p \bigl\|y_{B} \bigr\|^p \leq 1 + \tau^p.
\end{equation*}
Besides, using the fact that $(u+v)^\# = u^\# + v^\#$ for
disjointly supported elements $u,v \in L_p(\mu)$, it is clear
that $z_\theta^\# = \lambda^{- \frac{1}{q}} x_{A}^\# +
\overline{\theta} (1 - \lambda)^{-\frac{1}{q}} \tau^{p-1}
y_{B}^\#$. Using this and \eqref{eq:thm-nar-estimate} we can
write
\begin{align}\label{eq:thm-nar-essential}
(1+\tau^p)v(T)&\geq\max_{\theta\in\T}\left|\int_\Omega z_\theta^\# Tz_\theta\, d\mu\right|\\
&=\max_{\theta\in\T} \left| \lambda^{-1} \int_A x^\# T
x_{A} \,d \mu +  \theta\lambda^{- \frac{1}{q}} (1 -
\lambda)^{- \frac{1}{p}} \tau \int_A x^\# T y_{B} \,d \mu\right.\notag\\
&\qquad\quad \left.\phantom{\int_A}+ \overline{\theta} \lambda^{- \frac{1}{p}}(1 - \lambda)^{- \frac{1}{q}}
\tau^{p-1} \int_B y^\# T x_{A} \,d \mu + (1 -
\lambda)^{-1} \tau^p \int_B y^\# Ty_{B} \,d \mu \right|\notag\\
&\geq\max_{\theta\in\T}\left| \lambda^{-1} \int_A x^\# T
x_{A} \,d \mu + (1 -
\lambda)^{-1} \tau^p \int_B y^\# Ty_{B} \,d \mu+  \theta\lambda^{- \frac{1}{q}} (1 -
\lambda)^{- \frac{1}{p}} \tau \int_A x^\# T y_{B} \,d \mu\right.\notag\\
&\qquad\quad \left.\phantom{\int_A}+ \overline{\theta} \lambda^{\frac{1}{q}}(1 - \lambda)^{\frac{1}{p}}
\tau^{p-1}\|y\|^p \right|-\lambda^{- \frac{1}{p}}(1 - \lambda)^{- \frac{1}{q}}
\tau^{p-1}\Bigl| \int_B y^\# T x_{A} \,d \mu-\lambda(1-\lambda)\|y\|^p \Bigr|\notag\\
&\geq \max_{\theta\in\T}\left| \lambda^{-1} \int_A x^\# T
x_{A} \,d \mu + (1-\lambda)^{-1} \tau^p \int_B y^\# Ty_{B} \,d \mu+  \theta\lambda^{- \frac{1}{q}} (1-\lambda)^{- \frac{1}{p}} \tau \int_A x^\# T y_{B} \,d \mu\right.\notag\\
&\qquad\quad \left.\phantom{\int_A}+ \overline{\theta} \lambda^{\frac{1}{q}}(1 - \lambda)^{\frac{1}{p}}
\tau^{p-1}\|y\|^p \right|-\lambda^{- \frac{1}{p}}(1 - \lambda)^{- \frac{1}{q}}
\tau^{p-1}\eps\notag\\
&\geq\max_{\theta\in\T}\left|\theta\lambda^{- \frac{1}{q}} (1 -
\lambda)^{- \frac{1}{p}} \tau \int_A x^\# T y_{B} \,d \mu
+ \overline{\theta} \lambda^{\frac{1}{q}}(1 - \lambda)^{\frac{1}{p}}
\tau^{p-1}\|y\|^p \right| -
\tau^{p-1}\eps.\notag
\end{align}
From this point we study the real and the complex case
separately. For the complex case, we continue the estimation in
$\eqref{eq:thm-nar-essential}$ as follows
\begin{align*}
(1+\tau^p)v(T)&\geq\underset{\theta \in\T}{\max}\,\left|\theta\lambda^{- \frac{1}{q}} (1 -
\lambda)^{- \frac{1}{p}} \tau \int_A x^\# T y_{B} \,d \mu
+ \overline{\theta} \lambda^{\frac{1}{q}}(1 - \lambda)^{\frac{1}{p}}
\tau^{p-1}\|y\|^p \right| -\tau^{p-1}\eps\\
&=\left|\lambda^{- \frac{1}{q}} (1 -
\lambda)^{- \frac{1}{p}} \tau \int_A x^\# T y_{B} \,d \mu\right|\,
+ \,\left|\lambda^{\frac{1}{q}}(1 - \lambda)^{\frac{1}{p}}
\tau^{p-1}\|y\|^p \right| \, - \, \tau^{p-1}\eps \\
&\geq\lambda^{\frac{1}{q}}(1 - \lambda)^{\frac{1}{p}}
\tau^{p-1}\|y\|^p -\tau^{p-1}\eps\\
&\geq \lambda^{\frac{1}{q}}(1 - \lambda)^{\frac{1}{p}}
\tau^{p-1}(1-\eps)-\tau^{p-1}\eps.
\end{align*}
By the arbitrariness of $\eps$ we can write
$$
v(T)\geq\lambda^{\frac{1}{q}}(1 - \lambda)^{\frac{1}{p}}\frac{\tau^{p-1}}{1+\tau^p}
$$
for every $\tau>0$ and every $\lambda \in ]0,1[$ of the form
$\lambda = \frac{j}{2^n}$ where $j \in \{1, \ldots, 2^n-1\}$.
Since the diadic numbers are dense in $[0,1]$ and
$\underset{\lambda\in[0,1]}{\max}
\lambda^\frac{1}{q}(1-\lambda)^\frac{1}{p}=\kappa_p=\underset{\tau>0}{\max}\,\frac{\tau^{p-1}}{1+\tau^p}$\,,
the last inequality implies $v(T)\geq \kappa_p^2$ which finishes
the proof in the complex case.

In the real case, using (A) and (B) of Lemma~\ref{lem:nar3}, it
is easy to check that
\begin{equation*}
\lambda^{- \frac{1}{q}} (1 - \lambda)^{-
\frac{1}{p}} \tau \left| \int_A x^\# T y_{B} \,d \mu \right| \leq
\lambda^{- \frac{1}{q}} (1 - \lambda)^{- \frac{1}{p}} \tau \bigl\|
x^\#_{A} \bigr\|_q \bigl\| y_{B} \bigr\|_p \leq \lambda^{-
\frac{1}{q}} (1 - \lambda)^{- \frac{1}{p}} \tau \lambda^\frac{1}{q}
(1 - \lambda)^\frac{1}{p} = \tau
\end{equation*}
which, together with \eqref{eq:thm-nar-essential} and the choice of $y$, implies that
\begin{align*}
(1+\tau^p)v(T)&\geq\left|\lambda^{- \frac{1}{q}} (1 -
\lambda)^{- \frac{1}{p}} \tau \int_A x^\# T y_{B} \,d \mu
+\lambda^{\frac{1}{q}}(1 - \lambda)^{\frac{1}{p}}
\tau^{p-1}\|y\|^p \right| -
\tau^{p-1}\eps\\
&\geq\lambda^{\frac{1}{q}}(1 - \lambda)^{\frac{1}{p}}
\tau^{p-1}\|y\|^p -\lambda^{- \frac{1}{q}} (1 - \lambda)^{-
\frac{1}{p}} \tau \left| \int_A x^\# T y_{B} \,d \mu \right|-\tau^{p-1}\eps\\
&\geq \lambda^{\frac{1}{q}}(1-\lambda)^{\frac{1}{p}}
\tau^{p-1}(1-\eps)-\tau-\tau^{p-1}\eps.
\end{align*}
Hence, by the arbitrariness of $\eps$ we deduce that
$$
v(T)\geq\frac{\lambda^{\frac{1}{q}}(1-\lambda)^{\frac{1}{p}}
\tau^{p-1}-\tau}{1+\tau^p}\,
$$
for every $\tau>0$ and every $\lambda \in ]0,1[$ of the form
$\lambda = \frac{j}{2^n}$ where $j \in \{1, \ldots, 2^n-1\}$.
Taking supremum on $\lambda$, one has
$$
v(T) \geq \frac{\kappa_p \tau^{p-1} - \tau}{1 + \tau^p}
$$
for each $\tau > 0$, completing the proof.
\end{proof}

\section{Open problems}

\begin{problem}
{\slshape Calculate the numerical index of $L_p(\mu)$ for
$1<p<\infty$, $p\neq 2$.\,} As we commented in the introduction,
there are some estimations in the real case, but in the complex
case the knowledge about $n(L_p(\mu))$ is almost negligible. As a
conjecture, we think that $n(L_p(\mu))=M_p$ in the real case and
$n(L_p(\mu))=\kappa_p$ in the complex case.
\end{problem}

\begin{problem}
{\slshape Calculate the rank-one numerical index of
$L_p(\mu)$.\,} We conjecture that $n_1(L_p(\mu))=\kappa_p$ in
both the real and the complex cases (if the dimension of $L_p(\mu)$ is greater than $1$).
\end{problem}

\begin{problem}
{\slshape Is it true that the numerical index of $L_p(\mu)$
coincides with $n_{\rm nar}(L_p(\mu))$\,?\ } Let us comment
that for $Z=L_p([0,1],\ell_2)$ one has $n(Z)=n(\ell_2)<1$ and
$n_{\rm nar}(Z)=1$ since $Z$ has the so-called Daugavet
property. On the other hand, it is not difficult to show that
$n(\ell_p)$ coincides with the numerical index of compact
operators on $\ell_p$.
\end{problem}

\vspace{0.5cm}

\noindent\textbf{Acknowledgments:\ } The authors would like to
thank Rafael Pay\'{a} and Beata Randrianantoanina for fruitful
conversations concerning the matter of this paper.

\end{document}